\documentclass[rose]{article} 
\usepackage{amssymb,amsmath,latexsym, color}
 
\def\R {{\mathbb R}}

\def\N {{\mathbb N}}

\def\<{{\langle}}
\def\>{{\rangle}}

\def\proof{{\medskip\noindent {\bf Proof. }}}

\def\qed{{\hfill $\square$ \bigskip}}

\begin{document}

\newtheorem{thm}{Theorem}
\newtheorem{cor}[thm]{Corollary}
\newtheorem{prop}[thm]{Proposition}
\newtheorem{lem}[thm]{Lemma}
\newtheorem{remark}[thm]{Remark}
\newtheorem{Theo}[thm]{Definition}


\begin{center}
\LARGE{\bf Mixed Sub-Fractional  Brownian Motion}\\
\vspace{2mm}
Mounir Zili
 \end{center}                           %

\vspace{4mm}

\begin{quote}
{\bf Abstract:}
A new extension of the sub-fractional Brownian motion, and thus of the Brownian motion, is introduced. It is a linear combination of a finite number of sub-fractional Brownian motions, that we have chosen to call the mixed sub-fractional Brownian motion.
In this paper, we study some basic properties of this process, its non-Markovian and non-stationarity
characteristics,  the conditions under which it is a semimartingale, and the main features of its sample paths. We also show that this process could serve to get a good model of  certain phenomena,   taking not only  
 the sign (like in the case of the sub-fractional Brownian motion), but also  the strength  of dependence between the increments
of this phenomena into account. 
\end{quote}

Hausdorff dimension.\\


 \vspace{2mm}

\section{Introduction}

The sub-fractional Brownian motion (sfBm) is an extension of a Brownian motion (Bm), which was investigated in many papers (e.g.  \cite{TB}, \cite{Tu}).  It is a stochastic process 
$\displaystyle \xi^H = \{ \xi_t^H ; t \ge 0 \} $, defined on a  probability
space $\displaystyle (\Omega, F, {\mathbb P}) $ by:

\begin{equation}
 \forall t \in {\mathbb R}_+, \hspace{5mm}  \xi_t^H  
= \frac{B_t^H + B_{-t}^H}{\sqrt{2}}  ,
\label{eq:1}
\end{equation}

\noindent where $\displaystyle H \in ]0 , 1[$, and  $\displaystyle \{ B^H(t), t \in \R  \}$ is a fractional Brownian motion (fBm) on the whole real line; i.e. $B^H$ is a continuous and centered Gaussian process with covariance function 
\begin{equation}\label{eq:2} 
cov \Bigl( B^H(t), B^H(s) \Bigr) = 
1/2 \; \Bigl( \mid t \mid^{2H} + \mid s \mid^{2H} 
- \mid t - s \mid^{2H} \Bigr). 
\end{equation}
The index $H $ is called the Hurst parameter of $B^H$.  In some applications
(such as turbulence)  fBm is an adequate model for small increments,
but it seems to be inadequate for large increments (For more information about  fBm, see for example \cite{MaVa}). For this reason, $\xi^H$  may be an alternative to  fBm in some stochastic models. Moreover,
the sfBm arises from occupation time fluctuations of branching particle systems with Poisson initial condition \cite{TB}.

 In this paper, an extension of the sfBm  is introduced; it will be called the {\it mixed sub-fractional Brownian motion} (msfBm).  
 More precisely, for  $\displaystyle N \in \N \setminus\{ 0 \}$, $\displaystyle H= (H_1, H_2, ..., H_N) \in ]0,1[^N$ and 
$\displaystyle a= (a_1, a_2, ..., a_N) \in {\mathbb R} ^N\setminus \{ (0,..., 0)\}$, 
 the mixed sub-fractional Brownian motion  (msfBm) of parameters  $N$, $a$ and  $H$,
  is the  process
$\displaystyle S = \{ S_t^{H}(N,a) ; t \ge 0 \} =
  \{ S_t^{H}; t \ge 0 \} $, defined on the probability
space $\displaystyle (\Omega, F, {\Bbb P}) $ by:
\begin{equation}
\label{eq:3}
\forall t \in {\Bbb R}_+, \hspace{5mm}   S_t^{H}(N,a)  = \sum_{i=1}^N a_i \xi^{H_i} (t) 
\end{equation}

\noindent where $\displaystyle \Big( \xi^{H_i} \Big)_{i \in \{ 1, ..., N \} } $ is a family of independent sub-fractional Brownian motions of Hurst parameters $H_i$ defined on $\displaystyle (\Omega, F, {\Bbb P})$.

If $N= 1$,  and $\displaystyle a_1 = 1$,
\hspace{1mm}  $\displaystyle  S^{H}  = \xi ^H $ is a sub-fractional Brownian motion, 
and if  $N= 1$, $H_1 = 1/2$ and $\displaystyle a_1 = 1$, \hspace{1mm}  $\displaystyle S^{H}  $ is a standard  Brownian motion. 
So, the msfBm is more general, and this is a first reason which returns this process interesting to  be investigated.

In \cite{MZCE}, the authors studied on one hand some key properties of
the particular process $\displaystyle S^{(1/2, H_2)}(2,a)$ with $\displaystyle H_2 \in ]0,1[$ and $\displaystyle a= (a_1, a_2) \in {\mathbb R} ^2\setminus \{ (0,0)\}$, and on the other hand its martingale properties.

This paper is concerned with the study of the msfBm in general. Our first objective is to  extend the properties obtained in \cite{MZCE} to the general case. Particularly we  show that by a suitable choice of the parameters $a$ and $H$, the process $S^H$ could serve to get a good model of  certain phenomena,   taking not only  
 the sign (as in the case of the fBm and the sfBm), but also  the strength  of dependence between the increments
of the phenomena into account; and this is another main motivation to investigate this process.

Our second aim is to study some analytic and geometric   properties of the sample paths of $S^H$. Especially, we   investigate the Hausdorff dimensions of its graph, range and level sets.

The rest of this paper is organised as follows. Section $2$ contains basic properties of the msfBm, in particular the non-stationnarity,  the mixed-self-similarity \footnote{The {\it mixed-self-similarity} property was introduced by M. Zili in \cite{MZ}} and the non-Markovian properties.

In section $3$, we investigate some msfBm increments properties. We especially
 explicit   the covariances of the msfBm increments on non-overlapping intervals, and we analyze them to understand "how far" the msfBm  is from  a process with stationary increments. 
 
 Section $4$ is concerned with analytic and geometric properties of the sample paths of $S^H$. We show the H\"older-continuity and the non-differentiability of the tajectories, and we give the Hausdorff dimensions of 
 the graph, range and level sets. The methods used in our proofs are inspired by the papers \cite{Xiao}, \cite{Xiao2}
 and \cite{Ayache}, and are based particularly on the Frostman's theorem and the capacity argument (see \cite{Kah}).
 
In the last section, we investigate the semimartingale property, according to the values of the parameters 
$H$ and $a$. More precisely,
we show that the process $S^H(N,a)$ is a semimartingale, if and only if
 it exists $\displaystyle k_0 \in \{ 1, .., N \}; H_{k_0} = 1/2$, $\displaystyle a_{k_0} \neq 0$ and for every 
$\displaystyle i \in \{ 1,..., N \} \setminus \{ k_0 \} $, $\displaystyle H_i  \in \Big\{ 1/2 \Big\} \cup \Big]3/4, 1 \Big[$.

\section{The main properties}

Let us first recall some properties of the sfBm (see \cite{TB} and \cite{Tu} for proofs and for further information).

\begin{lem} \label{L:1}
The sfBm $\; (\xi_t^H)_{t \in {\mathbb R}_+}$ satisfies the following
properties:
\begin{itemize}
\item $\xi^H$ is a centered Gaussian process.
\item $ \displaystyle \forall s \in {{\mathbb R} }_+,  \forall t \in {{\mathbb R} }_+,$
\begin{equation} 
\label{eq:4} 
Cov \Big( \xi_t^H, \xi_s^H \Big) ) = s^{2H} + t^{2H} - 1/2  \Big( (s+t)^{2H} +  \mid t-s \mid ^{2H} \Big).\\
\end{equation}
\item For all $\displaystyle (s,t) \in {\Bbb R}_+^2; s \le t$,
\begin{equation}
\label{eq:5}
E \Big( \xi_t^{H} - \xi_s^{H} \Big) ^2  =   - 2^{2H-1} (t^{2H} + s^{2H}) +  (t+s)^{2H} +  (t -s) ^{2H}  .
\end{equation}
\item The increments of the smfBm are not stationary.
\end{itemize}
\end{lem}

As a first obvious consequence of Lemma \ref{L:1} we state:

\begin{lem} \label{L:2} 
{\it The msfBm  satisfies the following
properties: }
\begin{enumerate}
\item {\it $S_t^{H}$ is a centered Gaussian process.}
\item $\displaystyle \forall s \in {\Bbb R}_+$,
 $\displaystyle \forall t \in {\Bbb R}_+$, \\
\begin{equation}
\label{eq:6}
Cov \Big( S_t^{H}(a), S_s^{H}(a) \Big) =   \sum_{i=1}^N a_i^2 \Bigg[ t^{2H_i} + s^{2H_i} - 1/2 \Big[ ( s+t )^{2H_i} + \mid t -s \mid ^{2H_i} \Big]
\Bigg] .
\end{equation}
\item 
\begin{equation}
\label{eq:7}
 \forall t \in {\Bbb R}_+, 
\hspace{3mm} {\Bbb E}\Big( \big( S_t^H(a) \big) ^2 \Big) =  \sum_{i=1}^N a_i^2 \Big[ (2 - 2 ^{2H_i-1})  t^{2H_i} \Big].
\end{equation}

\end{enumerate}

\end{lem}

Let us now study the {\it mixed-self-similarity} property (see \cite{MZ}) of  the msfBm.

\begin{lem} \label{L:3}
 For any $h > 0$,  the processes 
 $\displaystyle  \big\{ S_{ht}^{H}(a) \big\}$ and   $\displaystyle   \Big\{
 S_t^{H} \Big( a_1 h^{H_1} ,   a_2 h^{H_2}, ..., a_Nh^{H_N} \Big) \Big\}$  have the same law.
\end{lem}

\proof  It is due to the fact that for fixed $h >0$, the processes 
$\displaystyle \{ S_{ht}^{H}(a) \} $ and
$ \displaystyle \Big\{ S_t^{H} \Big(  a_1 h^{H_1} ,   a_2 h^{H_2}, ..., a_Nh^{H_N} \Big) \Big\} $
  are Gaussian, centered and have the  same 
covariance function. 
\qed

The following lemma deals with the non-Markovian property of the msfBm. 

\begin{lem} \label{L:4} 
{\it For every  $\displaystyle a=(a_1, ..., a_N) \in {\Bbb R}^N $ and $\displaystyle H= (H_1, ..., H_N) \in \Big] 0; 1 \Big[^N $ such that  it exists $\displaystyle j \in \{ 1, ..., N \}; a_j \neq 0$ and $\displaystyle H_j \neq 1/2$, 
$\displaystyle (S_t^H(a))_{t \in {\Bbb R}_+}$ is not a Markovian process. } 

\end{lem}

\proof The process $S^H$ is a centered Gaussian. Moreover, it exists $\displaystyle j \in \{ 1, ..., N \}; a_j \neq 0$ and $\displaystyle H_j \neq 1/2$. Therefore, for all $t > 0$,
$$Cov \big( S^H_t, S^H_t \big)  = \sum_{i=1}^N  a_i^2\Big( (2-2^{2H_i-1}) \Big) t^{2H_i} > 0.$$
If $S^H$ were a Markovian process,  according
to \cite{Re},
for all $s < t < u$ we would have:
\begin{equation} 
\label{eq:8}
 Cov \Big( S^H_s,S^H_u \Big) Cov \Big( S^H_t, S^H_t \Big) = 
Cov \Big( S^H_s,S^H_t \Big) Cov \Big( S^H_t, S^H_u \Big) .
\end{equation}

Let us consider 
$$  H_{i_0} = \min \Big\{ H_i; i\in \{ 1,..., N \} ; a_i \neq 0 \Big\}  \hspace{2mm} {\rm and } \hspace{2mm}  H_{i_1} = \max \Big\{ H_i; i\in \{ 1,..., N \} ; a_i \neq 0 \Big\}  . $$

{\it First case:} If $H_j > 1/2$. In this case,  $\displaystyle H_{i_1} > 1/2 .$
By expression (\ref{eq:6}) and  equation (\ref{eq:8}) with  
$$1 < s = \sqrt{t} < t < u= t^2,$$
 we would have, 
\begin{equation} 
\label{eq:9}
\begin{array}{rcl}
&&\displaystyle \sum_{i=1}^N a_i^2 (2 - 2^{2H_i-1}) t^{2H_i} \\
\noalign{\vskip 2mm}
&& \displaystyle \times \sum_{i=1}^N a_i^2 \Bigg( t^{4H_i} + t^{H_i} - 1/2t^{4H_i}   \Big[ (1+t^{-3/2})^{2H_i} + (1- t^{-3/2})^{2H_i} \Big] \Bigg) \\
\noalign{\vskip 2mm}
& = & \displaystyle \sum_{i=1}^N a_i^2 \Bigg( t^{H_i} + t^{2H_i} - 1/2 t^{2H_i} \Big[ (1+t^{-1/2})^{2H_i} + (1-t^{-1/2})^{2H_i} \Big] \Bigg) \\
\noalign{\vskip 2mm}
&  & \displaystyle \times \sum_{i=1}^N a_i^2 \Bigg( t^{2H_i} + t^{4H_i} - 1/2 t^{4H_i} \Big[ (1+t^{-1})^{2H_i} + (1-t^{-1})^{2H_i} \Big] \Bigg) ,
\end{array}
\end{equation}

Since, when $h \rightarrow 0 $ we have
\begin{equation}
\label{eq:10}
 (1+h)^{2H_i} + (1-h)^{2H_i} = 2 + 2H_i(2H_i -1) h^2 + o(h^2),
 \end{equation}
 equations (\ref{eq:9}) and (\ref{eq:10}) imply that for $\displaystyle t \rightarrow \infty $,

\begin{equation} 
\label{eq:11}
\begin{array}{rcl}
&&\displaystyle  \sum_{i=1}^N a_i^2 (2 - 2^{2H_i-1}) t^{2H_i} \displaystyle \times \sum_{i=1}^N a_i^2 \Bigg(  t^{H_i} - H_i(2H_i-1)t^{4H_i-3}  + o\Big( t^{4H_i-3} \Big) \Bigg)  \\
\noalign{\vskip 2mm}
& - & \displaystyle \Bigg[ \sum_{i=1}^N a_i^2 \Bigg( t^{H_i}  - H_i(2H_i-1)t^{2H_i-1}  + o\Big( t^{2H_i-1} \Big) \Bigg) \\
\noalign{\vskip 2mm}
&  & \displaystyle \times \sum_{i=1}^N a_i^2 \Bigg( t^{2H_i}  - H_i(2H_i-1)t^{4H_i-2}  + o\Big( t^{4H_i-2} \Big)  \Bigg) \Bigg] \\
\noalign{\vskip 1mm}
& = & \displaystyle 0.
\end{array}
\end{equation}

The left member of equation (\ref{eq:11}) would tend to zero as t goes to infty; consequently, 
we would have
\begin{equation} 
\label{eq:12}
\lim_{t \rightarrow \infty } a_{i_1}^4 (1 - 2^{2H_{i_1}-1}) t^{3H_{i_1}}  = 0,
\end{equation}
which is true if, and only if  $H_{i_1} = 1/2$.    So $S^H$ is not a Markovian process. \\

{\it Second case:} If $H_j < 1/2$. In this case, 
$\displaystyle   H_{i_0} < 1/2 .$ Using the same technique as explained in the first case, but with 
 $$0< s = t^2 < t < u = \sqrt{t} < 1 \hspace{2mm} {\rm and} \hspace{2mm} t \rightarrow 0,$$
 we get
 \begin{equation} 
 \label{eq:13}
\lim_{t \rightarrow 0 } a_{i_0}^4 (1 - 2^{2H_{i_0}-1}) t^{6H_{i_0}}  = 0,
\end{equation}
which is true if, and only if  $H_{i_0} = 1/2$.  So  in this case too, $S^H$ is not a Markovian process. 
\qed

\section{Study of the msfBm increments}

In the following propostion, we will characterize the increments second moment of  the msfBm.

\begin{prop} \label{p:5} For all $\displaystyle (s,t) \in {\Bbb R}_+^2; s \le t$,

\begin{enumerate}
\item 
\begin{equation}
\label{eq:14}
 \begin{array}{rcl}
 && \displaystyle  E \Big( S_t^{H}(a) - S_s^{H}(a) \Big) ^2  \\
 \noalign{\vskip 3mm}
 &=& \displaystyle  \sum_{i=1}^N a_i^2 \Bigg( - 2^{2H_i-1} (t^{2H_i} + s^{2H_i}) +  (t+s)^{2H_i} +  (t -s) ^{2H_i} \Bigg) .
 \end{array}
   \end{equation}
 
\item 
\begin{equation}
\label{eq:15}
  \sum_{i=1}^N  a_i ^2 \gamma _i  (t-s)^{2H_i}   \le  E \Big( S_t^H(a) - S_s^H(a) \Big) ^2 \le \sum_{i=1}^N  a_i^2 \nu_i  (t-s)^{2H_i}  
  \end{equation}
where
\begin{equation}
\label{eq:16}
\gamma _i = \left\{ \begin{array}{rcl}
\displaystyle  2- 2^{2H_i-1} & if & \displaystyle  H_i > 1/2 \\
\noalign{\vskip 2mm}
\displaystyle  1 & if & \displaystyle H_i \le 1/2  \\
\end{array} 
\right.
\end{equation}
and
\begin{equation}
\label{eq:17}
\nu_i = \left\{ \begin{array}{rcl}
\displaystyle 1   & if & \displaystyle  H_i > 1/2 \\
\noalign{\vskip 2mm}
\displaystyle 2- 2^{2H-1} & if & \displaystyle H_i \le  1/2  .
\end{array} 
\right. 
\end{equation}

\end{enumerate}
 
 \end{prop}

\proof   The first result is due to equation (\ref{eq:5}) and to the fact that the processes $\xi^{H_i}$ are independent.   By the same fact, the second assertion is a direct consequence of  part $(3)$ of the Theorem on p. $407$ of \cite{TB}. 
\qed \\

{\bf Remark:} From proposition \ref{p:5}, it is clear that the msfBm does not have stationary increments, but this property is replaced by inequalities (\ref{eq:15}). \\

In the following Lemma, by an easy calculus and by equation (\ref{eq:6}), we  explicit   the covariances of  the msfBm increments on non-overlapping intervals. 

\begin{lem} \label{L:6}
 
If for $0 \le u < v \le s < t,$ we denote by 
$$ C_{u,v,s,t}= C_{u,v,s,t} (a) = Cov \Big( S_v^H(a)-S_u^H(a), S_t^H(a)-S_s^H(a) \Big) $$
then,
  \begin{equation}
\label{eq:19}
\begin{array}{rcl}
\displaystyle C_{u,v,s,t} &=& \displaystyle  \sum_{i=1}^N \frac{a_i^2}{2} \Big[ (t+u)^{2H_i} + (t-u)^{2H_i} + (s+v)^{2H_i} + (s-v)^{2H_i} \\
\noalign{\vskip 2mm}
&-& \displaystyle  (t+v)^{2H_i} - (t-v)^{2H_i} - (s+u)^{2H_i} - (s-u)^{2H_i} \Big] .
\end{array} 
\end{equation}
\end{lem}

As a first consequence of Lemma \ref{L:6},  we can  specify the sign of correlation between the increments of the msfBm, according to the values of $H$.

\begin{cor} \label{C:7}
For $0 \le u < v \le s < t,$
 \begin{enumerate}
\item $\displaystyle C_{u,v,s,t} = 0 $, if for every $\displaystyle i \in \{ 1, ..., N \}, H_i = 1/2 $,
 \item $\displaystyle C_{u,v,s,t} > 0$, if  for every $\displaystyle i \in \{ 1, ..., N \}, H_i > 1/2 $,
 \item $\displaystyle C_{u,v,s,t} < 0 $, if for every $\displaystyle i \in \{ 1, ..., N \}, H_i < 1/2 $.
 \end{enumerate}
\end{cor}

\proof The first assertion is trivial. To  check the second one,
let us  write
$$C_{u,v,s,t} = \sum_{i=1}^N \frac{a_i^2}{2} (g_i(t) - g_i(s))$$
where
$$g_i(t) = -(t+v)^{2H_i} - (t-v)^{2H_i} + (t+u)^{2H_i} + (t-u)^{2H_i}.$$
For every $\displaystyle i \in \{ 1, ..., N \},$ the function $g_i$ is differentiable and for every $t > 0$, 
$$g_i'(t) = 2H_i \Big( -(t+v)^{2H_i-1} - (t-v)^{2H_i-1} + (t+u)^{2H_i-1} + (t-u)^{2H_i-1} \Big) .$$
Since  $\displaystyle H_i > 1/2 $,  the function 
$\displaystyle x \longmapsto x^{2H_i-1}$ is concave and  $\displaystyle g_i'(t) > 0$.
Consequently, $\displaystyle g_i$ increases  and the second assertion holds.\\

The third result proof is similar.

\qed

As a second consequence of Lemma \ref{L:6}, it is easy to get the following main result:

\begin{cor} \label{C:8}
For $0 \le u < v \le s < t,$ $\displaystyle i \in \{ 1, ..., N \}$ and $\displaystyle (a_j)_{j \in \{ 1, ..., N \} \setminus \{ i \} }$, 
if $b$ and $c$ are two real numbers such that; $\displaystyle \mid b \mid \le \mid c \mid$, then
\begin{equation}
\label{eq:20}
C_{u,v,s,t}(a_1,...,a_{i-1},b,a_{i+1},..., a_N) <  C_{u,v,s,t}(a_1,...,a_{i-1},c,a_{i+1},..., a_N)   
\end{equation}
$\displaystyle \Big[ respectively \hspace{1mm} = ,\hspace{1mm}  >  \Big]$
if  \hspace{1mm} $H_i < 1/2$ \Big[ respectively \hspace{1mm} $H_i = 1/2$, \hspace{1mm} $H_i >  1/2$ \Big] .
\end{cor}

So, from Corollary \ref{C:8} we see that, for $0 \le u < v \le s < t,$ $\displaystyle i \in \{ 1, ..., N \}$ and $\displaystyle (a_j)_{j \in \{ 1, ..., N \} \setminus \{ i \} }$, if $\displaystyle H_i > 1/2$ $\Big[ $ respectively 
$\displaystyle H_i < 1/2$ $\Big]$,
\begin{itemize}
\item  the smaller $\Big[$  larger $\Big]$  $\mid b \mid $ is, the stronger
the dependence between the increments of  $\displaystyle S^H((a_1,...,a_{i-1},b,a_{i+1},..., a_N)$
is, 
\item  the larger $\Big[$  smaller $\Big]$ $\mid b \mid $ is, the weaker
 the dependence between the increments of  $\displaystyle S^H((a_1,...,a_{i-1},b,a_{i+1},..., a_N)$ is.
\end{itemize}

{\it Consequence:}  
In the modelling of a certain phenomena,  
we can choose $H=(H_1, ..., H_N)$ and $a=(a_1, ..., a_N)$
suitably in a such manner that   $\displaystyle \{ S_t^H(a) \} $
 permits to obtain a good model, taking not only  
 the sign (as in the case of fBm and sfBm), but also  the strength  of dependence between the increments
of the phenomena into account.

 In the following Lemma, we  prove that  the increments of  the msfBm are short-range dependent. For this, let us denote, for every nonnegative  real $x$ and nonnegative integer $n$, 
\begin{equation}
 \label{eq:21}
C(x,n) = C_{x,x+1,x+n,x+n+1} = Cov\Big( S^H_{x+1}-S^H_{x}, S^H_{x+n+1}-S^H_{x+n} \Big) .
\end{equation}

By Lemma \ref{L:6}, a Taylor expansion and an easy calculus we get:

\begin{lem} 
 \label{L:9} 
\begin{enumerate}
\item For every nonnegative integers $p$ and $n$ such that $n \ge 1$, we have
\[ \begin{array}{rcl} 
\displaystyle C(p, n) &=& \displaystyle  \sum_{i=1}^{N} \frac{a_i^{2}}{2} \Bigg[  (n+1)^{2H_i} - 2 n^{2H_i} + (n-1)^{2H_i} \\
\noalign{\vskip 2mm}
&& \displaystyle - (2p+n+2)^{2H_i} + 2 (2p+n+1)^{2H_i} - (2p+n)^{2H_i}\Bigg] .
\end{array} \]
\item For every $\displaystyle H \in ]0,1[^N $ and $\displaystyle p \in \N$, it holds that as $\displaystyle n \rightarrow \infty $, 
$$C(p,n) = \sum_{i=1}^N \Big[ 2(1-H_i)H_i(2H_i-1)(2p+1)a_i^2  n^{2H_i-3} + o(n^{2H_i-3}) \Big].$$
\item For every $\displaystyle H \in ]0, 1[^N $ and every $p \in \N $,
$$\sum_{n \ge 0} C(p, n) < \infty .$$
\end{enumerate}

\end{lem}

We can now analyze the function $\displaystyle x \longmapsto C(x,n)$ to understand "how far"  the msfBm is from  a process with stationary increments. Let us first consider,  for  $\displaystyle N \in \N \setminus\{ 0 \}$, $\displaystyle H= (H_1, H_2, ..., H_N) \in ]0,1[^N$ and 
$\displaystyle a= (a_1, a_2, ..., a_N) \in \R ^N\setminus \{ (0,..., 0)\}$, the process,  
$\displaystyle M^H = \{ M_t^{H}(N,a) ; t \ge 0 \}  $, defined on the probability
space $\displaystyle (\Omega, F, {\Bbb P}) $ by:
\begin{equation}
 \forall t \in {\Bbb R}_+, \hspace{5mm}  M_t^{H}(N,a) = \sum_{i=1}^N a_i B_t^{H_i}  
\label{eq:22}
\end{equation}
where the $B_t^{H_i}$'s are independent fractional Brownian motion with Hurst parameters $H_i$ defined on   
 $\displaystyle (\Omega, F, {\Bbb P}) $.\\
  This process is called  the mixed fractional Brownian motion (mfBm). It has been introduced  in $2001$ by Cheridito \cite{Cherid}, studied in $2007$ by M.Zili \cite{MZ} in the particular case where $N= 2$ and $H_1= 1/2$,  then extended in $2008$, by  Y.Miao-W.Ren-Z.Ren \cite{Miao} to the case where 
 $N= 2$ and $H_1 \in ]0,1[$. In $2009$, further remarks on the generalized form (\ref{eq:22}) of $M^H$ have been given by
 C. Th\"ale \cite{Tha}.  

The following Lemma is  due to \cite{MZ}:

\begin{lem} 
\label{L:10}
If for every positive real $x$, we denote by
\begin{equation} 
\label{eq:23}
R(x,n) = R_{x,x+1,x+n,x+n+1} = Cov\Big( M^H_{x+1}-M^H_{x}, M^H_{x+n+1}-M^H_{x+n} \Big), 
\end{equation}
 then
\begin{equation}
\label{eq:24}
 R(x,n ) = R(0,n) = \sum_{i=1}^N \frac{a_i^{2}}{2} \Bigg[  (n+1)^{2H_i} - 2 n^{2H_i} + (n-1)^{2H_i} \Bigg] .
 \end{equation}
\end{lem}

From Lemma \ref{L:10} it is obvious that 
the increments of the mfBm are stationnary. So, to meet our target, it suffices to
analyze the function $\displaystyle x \longmapsto C(x,n)$ to understand "how far"  the msfBm is from  the mfBm.

\begin{lem} \label{L:11}
For every $n$ it holds that as $x \rightarrow \infty $,
\begin{equation}
\label{eq:25}
C(x,n) = R(0,n)  - \sum_{i=1}^N a_i^2 2^{2H_i-1} H_i(2H_i-1) x^{2(H_i-1)} + o(x^{2(H_i-1)}),
\end{equation}
therefore $\displaystyle \lim_{x \rightarrow \infty} C(x,n) = R(0,n)$ for every $n$.
\end{lem}

\proof By Lemma \ref{L:9}, Lemma  \ref{L:10}  and  a Taylor expansion, we get the result.
\qed 

\section{Study of  the msfBm sample paths}

 In the sequel of the paper, we denote  by $i_0$ an integer such that
\begin{equation}
\label{eq:26}
i_0 \in \{ 1, ..., N \} \hspace{1mm} and \hspace{1mm}  H_{i_0} = \min \Big\{ H_i; i \in \{ 1, ..., N \} \hspace{1mm} {\rm and} \hspace{1mm} a_i \neq 0 \Big\} .
\end{equation}

In the following lemma, we check the continuity of the msfBm sample paths, and we even show that the parameter $H_{i_0}$, and consequently $H$, controls their regularity. \\

\begin{lem} 
 \label{L:12}  For all $T >0$ and  $ \displaystyle \gamma < H_{i_0}$,  the msfBm has a modification which sample-paths have 
a H\"{o}lder-continuity, with order $ \displaystyle \gamma $,
on the interval $\displaystyle [0; T]$.

\end{lem}

{\it Proof.} By  \cite{TB}, for every $\displaystyle i \in \{ 1, ..., N \}$, if $H_i \in ]0, 1[$ and $\displaystyle 0 < \epsilon < H_i $,  there exists a nonnegative random variable $\displaystyle G_{H_i, \epsilon , T}$  such that 
$\displaystyle E(G_{H_i, \epsilon , T}^p) < \infty $  for every $ p \ge 1$, 
and
$$\mid \xi ^{H_i}(t)- \xi ^{H_i}(s) \mid \le G_{H_i, \epsilon , T} \mid t-s \mid ^{H_i - \epsilon } \hspace{2mm} a.s,$$
for all $\displaystyle s,t \in [0,T]$.\\

So, for every $0 < \epsilon < H_{i_0} $, and $\displaystyle s,t \in [0; T]$
$$ \mid S_t^H - S_s^H \mid  = \Big| \sum_{i=1}^N  a_i  (\xi ^{H_i}(t)-\xi ^{H_i}(s)) \Big| \le G_{\epsilon , T} \mid t-s \mid ^{H_{i_0}- \epsilon } \hspace{2mm} a.s. $$
where
$$ G_{\epsilon , T} = \sum_{i=1}^N  \mid a_i \mid G_{H_i, \epsilon , T} T^{H_i-H_{i_0}}$$
for which we have clearly, $\displaystyle E(G_{\epsilon , T}^p) < \infty $  for every $ p \ge 1$. And the proof is finished.
\qed

In the next Lemma, we prove that with probability one, the msfBm sample paths  are not differentiable.

\begin{lem}  
\label{L:13} For every $\displaystyle H \in ] 0; 1[^N$, 
\begin{equation}
\label{eq:27}
\lim_{\epsilon \rightarrow 0^+} \sup_{t \in [t_0-\epsilon , t_0+ \epsilon ]} \Bigg| \frac{ S^H(t) - S^H(t_0) }{t-t_0}
\Bigg|  = + \infty ,
\end{equation}
with probability one for every $\displaystyle t_0 \in {\mathbb R} $.
 
\end{lem}

{\proof.} If for $m \in \N$, we denote $\displaystyle A^{(m)} = \cap_{n=1}^{+ \infty }A_n^{(m)} $ where, for $n \in \N$,
$$A_n^{(m)} = \Bigg\{ w \in \Omega; \sup_{t \in [t_0-\frac{1}{n} , t_0+\frac{1}{n} ]} 
\Bigg| \frac{ S^H(t) - S^H(t_0) }{t-t_0}
\Bigg|  > m \Bigg\} ,$$ 
to obtain the result it is enough to prove that $\displaystyle P \Big( \cap_{m=1}^{+ \infty } A^{(m)} \Big) = 1.$\\

On one hand,   the fact that each of the sequences $\displaystyle \Big( A_n^{(m)} \Big) _n$ and $\displaystyle \Big( A^{(m)} \Big) _m$
decreases, allows to write
$$P \Big( \cap_{m=1}^{+ \infty } A^{(m)} \Big)  = \lim_{m \rightarrow + \infty } P \Big(  A^{(m)} \Big) 
= \lim_{m \rightarrow + \infty } \lim_{n \rightarrow + \infty } P \Big(  A_n^{(m)} \Big) .$$
On another hand,
$$P (A_n^{(m)}) \ge P \Bigg( \Big| S^H(t_0 + \frac{1}{n} ) - S^H(t_0) \Big| > \frac{m}{n}  \Bigg) .$$
So to prove the lemma, it is enough to show that 
$$\forall m \in \N, \lim_{n \rightarrow + \infty }  P \Bigg( \Big| S^H(t_0 + \frac{1}{n} ) - S^H(t_0) \Big| \le \frac{m}{n}  \Bigg) = 0.$$
Since $\displaystyle S^H(t_0 + \frac{1}{n} ) - S^H(t_0)$ is a centered Gaussian random variable with variance
$$\sigma^2_n (t_0) = E \Big( S^H(t_0 + \frac{1}{n} ) - S^H(t_0) \Big)^2,$$
we have:
\[ \begin{array}{rcl}
 \displaystyle P \Bigg( \Big| S^H(t_0 + \frac{1}{n} ) - S^H(t_0) \Big| \le \frac{m}{n}  \Bigg) &=&\displaystyle \frac{1}{\sigma _n(t_0) \sqrt{2\pi }} \int_{-\frac{m}{n}}^{\frac{m}{n}} 
\exp \Big( - \frac{x^2}{2\sigma _n^2(t_0)} \Big) dx \\
\noalign{\vskip 3mm}
& \le & \displaystyle 2 \frac{m}{n} \times \frac{1}{\sigma _n(t_0) \sqrt{2\pi }}.

\end{array} \]

From  Proposition \ref{p:5},
$$\sigma^2_n(t_0) = \sum_{i=1}^N  a_i^2 \Bigg( - 2^{2H_i-1} ((t_0+ \frac{1}{n})^{2H_i} + t_0^{2H_i}) +  (2t_0+ \frac{1}{n})^{2H_i} +  \frac{1}{n^{2H_i}} \Bigg) , $$
then, for large $n$,  $\displaystyle \sigma _n^2(t_0) \approx \sum_{i=1}^N  a_i^2 n^{-2H_i} $.

So $\displaystyle \lim_{n \rightarrow + \infty } n^2 \sigma _n^2(t_0) =\lim_{n \rightarrow + \infty }
 \sum_{i=1}^N  a_i^2 n^{2-2H_i} = + \infty $, and consequently
 $$\lim_{n \rightarrow + \infty }  P \Bigg( \Big| S^H(t_0 + \frac{1}{n} ) - S^H(t_0) \Big| \le \frac{m}{n}  \Bigg) = 0.$$
 \qed

 Denoting the range of  the restriction of $S^H(a)$ on $[0,T]$ by
 \begin{equation}
 \label{eq:28}
 S^H([0,T]) = \Big\{ S^H(t); t \in [0,T] \Big\} ,
 \end{equation}
 its graph by
 \begin{equation}
 \label{eq:29}
 Grf_T S^H(a)  = \Big\{ (t, S^H_t(a)); t \in [0, T] \Big\},
 \end{equation}
  the graph of $S^H(a)$ by
   \begin{equation}
 \label{eq:30}
 Grf S^H(a)  = \Big\{ (t, S^H_t(a)); t \in [0, + \infty [  \Big\},
 \end{equation}
and the level set  of the restriction of  $S^H(a)$ on $[\epsilon , T]$ by
\begin{equation}
\label{eq:31}
L_x^\epsilon = \Big\{ t \in [\epsilon , T]; S^H(t) = x \Big\} ,
\end{equation}
where  $\displaystyle T > \epsilon > 0$,
the aim of the sequel of this section, is to study the  Hausdorff dimensions of the sets defined in (\ref{eq:28}), 
 (\ref{eq:29}), (\ref{eq:30}) and (\ref{eq:31}).
 
  Let us first recall briefly the definition of Hausdorff dimension. For each $\alpha > 0$, $E \subset {\Bbb R}^d$, the $\alpha -$ dimensional Hausdorff measure of $E$ is defined by
\begin{equation}
\label{eq:32}
\textsl{M}^\alpha (E) = \lim_{\delta \rightarrow 0} inf \Bigg\{ \mid E \mid^\alpha ; E \subset \cup_{k=1}^{\infty } 
E_k; \mid E_k \mid < \delta \Bigg\} ,
\end{equation}

where $\mid E_k \mid $  is the diameter of the set $E_k$ and the infinimum is taken over all coverings $(E_k)_{k \in {\Bbb N}}$ of $E$. The Hausdorff dimension of $E$ is defined by

\begin{equation}
\label{eq:33}
dim E = inf \{ \alpha > 0; \textsl{M}^\alpha (E) = 0 \} = sup \{ \alpha > 0; \textsl{M}^\alpha (E) = + \infty \} .
\end{equation}

\begin{lem}
\label{L:14}
The Hausdorff dimension of the graph of $S^H(a)$ equals $\displaystyle 2 - H_{i_0}$ with probability $1$, where $H_{i_0}$ is defined by (\ref{eq:26}).
\end{lem}

\proof  By Lemma \ref{L:12},  for all $T >0$, the msfBm has a modification which sample-paths have 
a H\"{o}lder-continuity, with order $ \displaystyle \gamma < H_{i_0}$ on the interval $\displaystyle [0; T]$. 
So, according to \cite{Xiao}, for all $T > 0$, with probability $1$,  
$$dim Grf_T S^H(a)  \le 2 - H_{i_0},$$
which implies that
  $$dim Grf S^H(a)  \le 2 - H_{i_0}.$$
  
  Now, thanks to the Frostman's Theorem (see e.g. \cite{FAL}) to obtain the second inequality we only need  to show that for every $T >0$,  the occupation measure $\nu$ of $\displaystyle t \longmapsto (t, S_t^H(a))$, when $t$ is restricted to the interval $\displaystyle [0;T]$, has, with probability $1$, a finite $u-$dimensional energy, for any $\displaystyle u \in \Big] 1, 2 - H_{i_0} \Big[ .$ More precisely, for any Borel set $\displaystyle A \subset {\Bbb R}^2$, $\nu (A)$ is defined as the integral
  \begin{equation}
  \label{eq:34}
  \nu (A) = \int_0^T {\bf 1}_{\{ (t,S_t^H(a)) \in A \} } dt, 
  \end{equation}
where, for every set $U \subset {\Bbb R}^2$, $\displaystyle {\bf 1}_{U}$ denotes the characteristic function of  the set $U$, and we need to prove that with probability $1$ the integral
\begin{equation}
\label{eq:35}
\int_{Grf_T S^H(a)} \int_{Grf_T S^H(a)} \mid x -y \mid ^{-u} \nu (dx) \nu (dy)
\end{equation}
is finite. By a monotone class argument this is easily seen to be equivalent to

\begin{equation}
\label{eq:36}
\int_0^T \int_0^T \Big( \mid s -t \mid + \mid S_s^H(a) - S_t^H(a) \mid \Big) ^{-u} ds dt < + \infty,
\end{equation}

which follows from

\begin{equation}
\label{eq:37}
\int_0^T \int_0^T {\Bbb E} \Bigg( \Big( \mid s -t \mid + \mid S_s^H(a) - S_t^H(a) \mid \Big) ^{-u} \Bigg) ds dt < + \infty .
\end{equation}

In order to get (\ref{eq:37}), we need the following preliminary lemma.

\begin{lem}
\label{L:15}
For all $\displaystyle (s,t) \in {\Bbb R} \times {\Bbb R}, s \neq t$ and for every real $u > 1$, we have
\begin{equation}
\label{eq:38}
{\Bbb E} \Bigg( \Big( \mid s-t \mid + \mid S_t^H(a) -S_s^H(a) \mid \Big) ^{-u} \Bigg) \le
c \mid t-s \mid ^{1-u} \sigma ^{-1} (s,t),
\end{equation}
where 
\begin{equation}
\label{eq:39}
\sigma ^2(s,t) = {\Bbb E}  \Big(S_t^H(a) -S_s^H(a) \Big) ^2 
\end{equation}
and $c > 0$ is a constant.
\end{lem}

\proof (of Lemma \ref{L:15}) We have

\[ \begin{array}{rcl}
&& \displaystyle {\Bbb E} \Bigg( \Big( \mid s-t \mid + \mid S_t^H(a) -S_s^H(a) \mid \Big) ^{-u} \Bigg) \\
\noalign{\vskip 2mm}
&=& \displaystyle \frac{1}{\sigma (s,t) \sqrt{2\pi}} \int_{\Bbb R} \Big( \mid t-s \mid + \mid x \mid \Big) ^{-u}
\exp \Bigg( - \frac{x^2}{2 \sigma^2(s,t)} \Bigg) dx \\
\noalign{\vskip 2mm}
&=& \displaystyle \frac{2}{\sigma (s,t) \sqrt{2\pi}} \int_{0}^{\mid t-s \mid } \Big( \mid t-s \mid + \mid x \mid \Big) ^{-u}
\exp \Bigg( - \frac{x^2}{2 \sigma^2(s,t)} \Bigg) dx \\
\noalign{\vskip 2mm}
&& \displaystyle +  \frac{2}{\sigma (s,t) \sqrt{2\pi}} \int_{\mid t-s \mid }^{+ \infty } \Big( \mid t-s \mid + \mid x \mid \Big) ^{-u}
\exp \Bigg( - \frac{x^2}{2 \sigma^2(s,t)} \Bigg) dx \\
\noalign{\vskip 2mm}
&\le & \displaystyle \frac{2}{\sigma (s,t) \sqrt{2\pi}} \Bigg( \int_{0}^{\mid t-s \mid }  \mid t-s \mid ^{-u} dx +
\int_{\mid t-s \mid }^{+ \infty }  \mid x \mid ^{-u} dx \Bigg) \\
\noalign{\vskip 2mm}
&\le & \displaystyle \frac{2}{\sigma (s,t) \sqrt{2\pi}} \Bigg( \mid t-s \mid ^{1-u} +
\frac{1}{u-1} \mid t-s \mid ^{1-u} \Bigg) \\
\noalign{\vskip 2mm}
&\le & \displaystyle 
c \mid t-s \mid ^{1-u} \sigma ^{-1} (s,t),
\end{array} \]
with $$c = \frac{\sqrt{2}}{\sqrt{\pi}} 
\frac{u}{u-1}  .$$
\qed

Let us now prove (\ref{eq:37}). By  Lemma \ref{L:15} then by Proposition \ref{p:5} we get:

\[ \begin{array}{rcl}
&& \displaystyle \int_0^T \int_0^T {\Bbb E} \Bigg( \Big( \mid s -t \mid + \mid S_s^H(a) - S_t^H(a) \mid \Big) ^{-u} \Bigg) ds dt \\
\noalign{\vskip 2mm}
& \le & \displaystyle \int_0^T \int_0^T c \mid t-s \mid ^{1-u} \sigma ^{-1} (s,t) ds dt \\
\noalign{\vskip 2mm}
& \le & \displaystyle \int_0^T \int_0^T \frac{c}{a_{i_0} \sqrt{\gamma _{i_0}}}  \mid t-s \mid ^{1-u-H_{i_0}}  ds dt ,
\end{array} \]
where $\gamma _{i_0}$ is defined by (\ref{eq:16}).
Since  $\displaystyle u \in \Big] 1, 2 - H_{i_0} \Big[$, it is easy to check that

\begin{equation}
\label{eq:40}
\int_0^T \int_0^T   \mid t-s \mid ^{1-u-H_{i_0}}  ds dt  < + \infty ,
\end{equation}

which achieves the proof.
\qed

\begin{lem}
\label{L:16}
The Hausdorff dimension of the range  $S^H([0,T])$ equals $\displaystyle 1$ with probability $1$. 
\end{lem}

\proof
We have clearly $\displaystyle dim S^H([0,T]) \le 1 \hspace{2mm} a.s.,$ so we only need to prove that
$$1 \le dim \hspace{1mm} S^H([0,T]) \hspace{2mm} a.s.$$

Note that for $\displaystyle \epsilon \in ]0, T[$,
$$ dim S^H([0,T]) \ge dim S^H([\epsilon ,T]),$$
and that for any standard normal variable $X$ and $\displaystyle 0 < \gamma < 1$, we have 
\begin{equation}
\label{eq:41}
E(\mid X \mid ^{- \gamma }) < \infty.
\end{equation}

Hence by Frostman's theorem (see e.g. \cite{FAL}), it is sufficient to show that for all $\displaystyle 0 < \gamma < 1$,
\begin{equation}
\label{eq:42}
E_\gamma = \int_{\epsilon }^{T} \int_{\epsilon }^{T} E \Big( \mid S^H(s)-S^H(t) \mid ^{-\gamma } \Big) ds dt < + \infty.
\end{equation}

From Proposition \ref{p:5}, we see that  there exist positive and finite constants $c_1$ and $c_2$ such that for all $\displaystyle s,t \in [0, T]$,
\begin{equation}
\label{eq:43}
c_1 \mid s-t\mid ^{2H_{i_0}} \le E\Big( (S^H(s)-S^H(t))^2 \Big) \le c_2 \mid s-t\mid ^{2H_{i_0}}.
\end{equation}

So by (\ref{eq:42}) and (\ref{eq:43}), it exists a positive and finite constant $c_3$ such that

\begin{equation}
\label{eq:44}
E_\gamma \le c_3 \int_{\epsilon }^{T} \int_{\epsilon }^{T}   \mid s-t\mid ^{-\gamma H_{i_0}}  ds dt.
\end{equation}

Since $\displaystyle 0 < \gamma H_{i_0} < 1,$ the second member of the inequality (\ref{eq:44}) is finite and we get the result.  
\qed

The following lemma is necessary for the study of the Hausdorff dimension of the level set $L_x^\epsilon $.

\begin{lem}
\label{L:17}
If we denote by $\displaystyle Var(Y \mid Z)$ the conditional variance of $Y$ given $Z$, there exists a constant $c > 0$ such that for all $s,t \in I$,
\begin{equation}
\label{eq:45}
Var(S^H(t) \mid S^H(s) ) \ge c \mid s-t \mid^{2H_{i_0}}.
\end{equation}

\end{lem}

\proof Since the conditional variance in (\ref{eq:45}) is the square of the $L^2({\Bbb P})-$ distance of $S^H(t)$ from the subspace generated by $S^H(s)$, we have
\begin{equation}
\label{eq:46}
Var(S^H(t) \mid S^H(s) ) = \inf_{b \in {\Bbb R}} {\Bbb E} \Big( S^H(t) - b S^H(s) \Big) ^2. 
\end{equation}

So, from the definition (\ref{eq:3}) of the msfBm, and from the independence of the sfBm's $\xi^{H_i}$ we can write

\begin{equation}
\label{eq:47}
\begin{array}{rcl}
Var(S^H(t) \mid S^H(s) ) & = &  \displaystyle \inf_{b \in   {\Bbb R} }
\sum_{i=1}^{N} {\Bbb E} \Big( \xi^{H_i}(t) - b \xi^{H_i}(s) \Big) ^2 \\
\noalign{\vskip 2mm}
& \ge & \displaystyle \inf_{b \in   {\Bbb R} } {\Bbb E} \Big( \xi^{H_{i_0}}(t) - b \xi^{H_{i_0}}(s) \Big) ^2 \\
\noalign{\vskip 2mm}
&= & \displaystyle Var(\xi^{H_{i_0}}(t) \mid \xi^{H_{i_0}}(s) ).
\end{array}
\end{equation}

Thanks to Yan and Shen [Theorem $2.1$, \cite{Yan}], we know that it exists a constant $c > 0$ such that
\begin{equation}
\label{eq:48}
Var(\xi^{H_{i_0}}(t) \mid \xi^{H_{i_0}}(s) ) \ge c \mid t-s \mid ^{2H_{i_0}}.
\end{equation}

Equations (\ref{eq:47}) and (\ref{eq:48}) complete the proof.\\

\begin{lem}
\label{L:18}
For every $\displaystyle x \in {\Bbb R}, $ and $\displaystyle 0 < \epsilon < T$, with positive probability 
\begin{equation}
\label{eq:49}
dim_H(L_x^\epsilon  ) = 1 - H_{i_0},
\end{equation}
where $H_{i_0}$ is defined by (\ref{eq:26}).

\end{lem}

\proof
For an integer $\displaystyle n \ge 1$, devide the interval $\displaystyle [ \epsilon , T]$ into  $m_n$ sub-intervals $I_{n,l}$ of length $\displaystyle n^{-1/H_{i_0}}$. Then 
\begin{equation}
\label{eq:50}
 m_n \le T \times n^{1/H_{i_0}}.
 \end{equation}
  Let $0 < \delta < 1$ be fixed and let $\tau_{n,l} = 
\epsilon + l n^{-1/H_{i_0}}$. By Lemma $3.3$ of \cite{WuXiao} and Lemma $2.2$ of \cite{Talag} we get:
\begin{equation}
\label{eq:51}
\begin{array}{rcl}
\displaystyle {\bf P} \Big\{ x \in S^H(I_{n,l}) \Big\} & \le & \displaystyle  {\bf P} \Big\{ \max_{s,t \in I_{n,l}} \mid S^H(s) - S^H(t) \mid \le
n^{-(1-\delta )}; x \in S^H(I_{n,l})  \Big\} \\
\noalign{\vskip 3mm}
&& \displaystyle  + {\bf P} \Big\{ \max_{s,t \in I_{n,l}} \mid S^H(s) - S^H(t) \mid >
n^{-(1-\delta )}  \Big\} \\
\noalign{\vskip 3mm}
&\le &  \displaystyle  {\bf P} \Big\{  \mid S^H(\tau_{n,l}) - x \mid \le
n^{-(1-\delta )} \Big\} + \exp \Big( - c_1 n^{2\delta } \Big) \\
\noalign{\vskip 3mm}
& \le & \displaystyle  c_2 n^{-(1-\delta)}.
\end{array} 
\end{equation}

 Define a covering $\displaystyle \{ I'_{n,l} \} $ of $\displaystyle L_x $ by $\displaystyle I'_{n,l} = I_{n,l}$ if $\displaystyle x \in S^H(I_{n,l})$ and 
 $\displaystyle I'_{n,l} = \emptyset $ otherwise. Denote the number of such sets $\displaystyle \{ I'_{n,l} \} $ by $M_n$. By (\ref{eq:50}) and (\ref{eq:51}) we get:\\
 
  \begin{equation}
 \label{eq:52}
 \begin{array}{rcl}
 \displaystyle E(M_n) &\le & \displaystyle  E \Big( m_n \times {\bf 1}_{ \{ x \in S^H(I_{n,l}) \} } \Big) \\
 \noalign{\vskip 3mm}
 & \le & \displaystyle  T \times n^{1/H_{i_0}} \times  {\Bbb P}  \Big\{ x \in S^H(I_{n,l}) \Big\}  \\
 \noalign{\vskip 3mm}
 & \le & \displaystyle  c_3 \times n^{1/H_{i_0}-1 + \delta }   
 \end{array}
 \end{equation}
 where $c_3$ is a positive constant.\\
 
 Let $\displaystyle \eta = 1/H_{i_0} - (1-2\delta ).$ We consider the sequence of integers $\displaystyle n_i = 2^i (i \ge 1)$.
 Then by (\ref{eq:52}), the Markov inequality and the Borel-Cantelli lemma we see that almost surely 
 $\displaystyle M_{n_i} \le c_4 n_i^{\eta }$ for all $i$ large enough. This implies that $\displaystyle dim_H L_x^\epsilon  \le H_{i_0} \eta $ almost surely. Letting $\delta \downarrow 0$ along rational numbers, we  get
 \begin{equation}
 \label{eq:53}
 dim_H L_x^\epsilon \le 1 - H_{i_0} \hspace{2mm} a.s.
 \end{equation}
 
 To prove the lower bound for $\displaystyle dim_H L_x^\epsilon $ in (\ref{eq:49}), we consider $\delta >0 $ a small constant such that
 \begin{equation}
 \label{eq:54}
 \gamma = 1 - H_{i_0} (1+ \delta ) > 0.
  \end{equation}
  
  Note that if we can prove that there is a constant $c_5 > 0$, independent of $\delta$, such that
  \begin{equation}
 \label{eq:55}
 {\bf P} \Big\{ dim_H L_x^\epsilon \ge \gamma \Big\} \ge c_5,
  \end{equation}
  then the lower bound in (\ref{eq:49}) will follow by letting $\delta \downarrow 0$.\\

 Our proof of (\ref{eq:55}) is based on the capacity argument due to Kahane [see \cite{Kah}]. Similar methods have been used in \cite{Adler},  \cite{Ayache}, \cite{Testard} and \cite{Xiao}.
 
 Let $\displaystyle {\cal M}_{\gamma }^+$ be the space of all non-negative measures on ${\Bbb R}$ with finite $\gamma -$ energy. It is known [cf. \cite{Adler}] that $\displaystyle {\cal M}_{\gamma }^+$ is a complete metric space under the metric
 \begin{equation}
 \label{eq:56}
 \parallel \mu \parallel _{\gamma } = \int_{\Bbb R} \int_{\Bbb R} \frac{\mu (dt) \mu (ds)}{\mid t-s \mid ^\gamma } .
 \end{equation}
 
 We define a sequence of random positive measures $\displaystyle \mu_n:= \mu_n (x, .)$ on the Borel sets $C$ of 
 $\displaystyle [\epsilon , T]$ by
 \begin{equation}
 \label{eq:57}
 \begin{array}{rcl}
 \mu_n(C)& =& \displaystyle \int_C \sqrt{2 \pi } \exp \Big( - (S^H(t) -x)^2/2 \Big) dt \\
 \noalign{\vskip 2mm}
 & =& \displaystyle \int_C \int_{\Bbb R}  \exp \Big( - \xi^2/2 + i \xi (S^H(t)-x) \Big) d\xi dt. \\
 \noalign{\vskip 2mm}
 \end{array}
 \end{equation}
 
 It follows from \cite{Kah} (p. $206$) or \cite{Testard} (p.$17$) that if there exist positive and finite constants $c_6, c_7$ and $c_8$ such that
 \begin{equation}
 \label{eq:58}
 {\Bbb E} (\parallel \mu_n \parallel ) \ge c_6, \hspace{3mm}  {\Bbb E} (\parallel \mu_n \parallel ^2) \le c_7,
 \end{equation}
 
 \begin{equation}
 \label{eq:59}
 {\Bbb E} (\parallel \mu_n \parallel _\gamma ) \le c_8,
 \end{equation}
 where $\displaystyle \parallel \mu_n \parallel = \mu_n \Big( [\epsilon , T ] \Big) $ denotes the total mass of $\mu_n$, then there is a subsequence of $\{ \mu_n \}$, say $\{ \mu_{n_k}  \}$, such that $\mu_{n_k} \rightarrow \mu$
 in $\displaystyle {\cal M}_{\gamma }^+$ and $\mu$ is strictly positive with probability $\ge c_6^2/(2c_7)$.
 In this case, it follows from (\ref{eq:57}) that $\mu$ has its support in $L_x^\epsilon $ almost surely. Moreover, (\ref{eq:59}) and the monotone convergence theorem together imply that the $\gamma-$energy of $\mu$ is finite. Hence Frostman's theorem yields (\ref{eq:55}) with $\displaystyle c_5= c_6^2/(2c_7)$.\\
 
 It remains to verify (\ref{eq:58}) and (\ref{eq:59}). By Fubini's theorem we have
 
 \begin{equation}
 \label{eq:60}
 \begin{array}{rcl}
 \displaystyle {\Bbb E} (\parallel \mu_n \parallel ) & = & \displaystyle \int_{\epsilon }^T \int_{\Bbb R}
 \exp(-i\xi x) \exp(-\xi^2/2) {\Bbb E} \Big( \exp (i\xi S^H(t)) \Big)d \xi dt \\
 \noalign{\vskip 2mm}
 &= & \displaystyle \int_{\epsilon }^T  \int_{\Bbb R}
 \exp(-i \xi x)   \exp \Big( -(1+ \sigma ^2(t))\xi^2 /2 \Big)  d \xi dt \\
 \noalign{\vskip 2mm}
 &= & \displaystyle \int_{\epsilon }^T \sqrt{2\pi / (1+ \sigma ^2(t))} \exp \Big( - x^2/(2(1+ \sigma^2(t))) \Big) dt \\
 \noalign{\vskip 2mm}
 &\ge & \displaystyle \int_{\epsilon }^T \sqrt{2\pi / (1+ \sigma ^2(t))} \exp \Big( - x^2/2 \sigma^2(t) \Big) dt 
 := c_6.
 \end{array}
 \end{equation}
 
 Denote by $I_2$ the identity matrix of order $2$ and $Cov(S^H(s), S^H(t))$ the covariance matrix of the random vector
 $(S^H(s), S^H(t))$. Let $\displaystyle \Gamma = I_2 + Cov (S^H(s), S^H(t))$ and $(\xi , \eta )'$ be the transpose of the row vector $(\xi , \eta )$. Then

 \begin{equation}
 \label{eq:61}
 \begin{array}{rcl}
 \displaystyle {\Bbb E} (\parallel \mu_n \parallel ^2 ) & = & \displaystyle \int_{\epsilon }^T \int_{\epsilon }^T \int_{\Bbb R} \int_{\Bbb R}
 \exp \Big( -i(\xi + \eta ) x \Big) \exp \Big( -(\xi , \eta ) \Gamma (\xi , \eta ) ' / 2 \Big) d \xi d \eta ds dt  \\
 \noalign{\vskip 2mm}
 &= & \displaystyle \int_{\epsilon }^T  \int_{\epsilon }^T 
 \Big( 2\pi / \sqrt{det \Gamma } \Big) \exp \Big( -  (x , x ) \Gamma ^{-1} (x , x ) ' / 2 \Big) ds dt \\
 \noalign{\vskip 2mm}
 &\le  & \displaystyle \int_{\epsilon }^T  \int_{\epsilon }^T 
  2\pi / \sqrt{det Cov(S^H(s), S^H(t))}   ds dt .
 \end{array}
 \end{equation}

 Since 
 \begin{equation}
 \label{eq:62}
  det Cov( S^H(s), S^H(t)) = Var (S^H(s)) Var \Big( S^H(t) \mid S^H(s) \Big) ,
  \end{equation}
  from conditions 
 (\ref{eq:7}) and (\ref{eq:45})  we get: for all $\displaystyle s,t \in [ \epsilon , T ],$
 
 \begin{equation}
 \label{eq:63}
 det Cov(S^H(s), S^H(t)) \ge c_9 \mid s-t\mid ^{2H_{i_0}} .
  \end{equation}
  where $c_9$ denotes a positive constant.
  Combining (\ref{eq:62}) and (\ref{eq:63}) we obtain:
  
  \begin{equation}
 \label{eq:64}
  {\Bbb E} (\parallel \mu_n \parallel ^2 )  \le 
  \displaystyle c_{10} \int_{\epsilon }^T  \int_{\epsilon }^T 
  2\pi / \mid s-t\mid ^{H_{i_0}}   ds dt .
 \end{equation}
 where $c_{10}$ denotes a positive constant.
 Since $\displaystyle H_{i_0} \in ]0 , 1 [$, the last integral is a finite constant, and consequently
 \begin{equation}
 \label{eq:65}
 {\Bbb E} (\parallel \mu_n \parallel ^2 ) \le c_7
 \end{equation}
 with 
 $\displaystyle c_7 = c_{10} \int_{\epsilon }^T  \int_{\epsilon }^T 
  2\pi / \mid s-t\mid ^{H_{i_0}}   ds dt. $ \\
 
 Similar to (\ref{eq:61}), we have
 \begin{equation}
 \label{eq:66}
 \begin{array}{rcl}
 && \displaystyle {\Bbb E} (\parallel \mu_n \parallel _\gamma  ) \\
 \noalign{\vskip 3mm}
 & = & \displaystyle \int_{\epsilon }^T \int_{\epsilon }^T \mid s-t \mid ^{-\gamma } ds dt \int_{\Bbb R} \int_{\Bbb R}
 \exp \Big( -i(\xi + \eta ) x \Big) \exp \Big( -(\xi , \eta ) \Gamma (\xi , \eta ) ' / 2 \Big) d \xi d \eta   \\
 \noalign{\vskip 3mm}
 &\le  & \displaystyle \displaystyle c_{11} \int_{\epsilon }^T  \int_{\epsilon }^T 
  2\pi \mid s-t\mid ^{-\gamma - H_{i_0}}   ds dt \\
 \end{array}
 \end{equation}
 where $c_{11}$ denotes a positive constant. 
 Since 
 $\displaystyle - \gamma - H_{i_0}  = -1 + H_{i_0} \delta \in ]-1, 0[,$ we get
 $\displaystyle {\Bbb E} (\parallel \mu_n \parallel _\gamma  ) < c_8$ with
 $$c_8 =  c_{11} \int_{\epsilon }^T  \int_{\epsilon }^T 
  2\pi \mid s-t\mid ^{-\gamma - H_{i_0}}   ds dt,$$
  which completes the proof
 \qed 
 
\section{Study of the semimartingale property}

In this section, we will discuss for which values of the Hurst parameter $H$, $S^H(N,a) $ is a semimartingale. Let us first specify that 
  in this paper,  for a stochastic process $\{X_t, 0 \le t  \le T \}$,  we denote by $\displaystyle \mathcal{F}^X =({\cal F}_{t}^X)_{0 \le t \le T}$ the own filtration of $X$, and we call $X$ a semimartingale if it is a semimartingale with respect to  $\displaystyle \mathcal{\overline{F}}^{X}$, the smallest filtration that contains $\displaystyle \mathcal{F}^X $ and satisfies the usual assumptions.
Let us start our study by the following lemma.

\begin{lem}  \label{L:20}
Let us denote, for every $t> 0$,  by $\displaystyle <S^H(N,a)>_t$ \; [respectively $\displaystyle V(S^H )_t$ ] \; the quadratic variation \; [respectively the variation]\;  of the process $\displaystyle S^H(N,a)$ on the interval $[0,t]$, and let us recall that $H_{i_0}$ denotes the parameter  defined by (\ref{eq:26}). For every $t > 0$,\\
 
\noindent 
$1.$ if  $ \displaystyle  H_{i_0} < 1/2$ then,  \;
$\displaystyle  < S^H(N,a) >_t 
\stackrel{a.s.}{=} + \infty ,$ \\

\noindent  $2.$ if,   $\displaystyle H_i  > 1/2 $ for every $\displaystyle i \in \{ 1, ..., N \}$,  then 
$$ < S^H(N,a) >_t \stackrel{a.s.}{=} 0 \; \hspace{2mm} and \; \hspace{2mm}  V( S^H )_t \stackrel{a.s.}{=}  \infty . $$
\end{lem}

\proof $1.$ For any $\displaystyle n \in {\mathbb N} \setminus \{ 0\}$, $ \displaystyle p \in {\mathbb N} \setminus \{ 0\}$  and  $\displaystyle t> 0$, we denote
\begin{equation}
\label{eq:504}
A_{n,p} =  \sum_{j=1}^n \; \Big|  S_{\frac{jt}{n}}^H(N,a) - S_{\frac{(j-1)t}{n}}^H(N,a)  \Big| ^p ,
\end{equation} 

Assume that $<S^H(N,a)>_t < \infty \; a.s.$ Then $A_{n,2} \rightarrow <S^H(N,a)>_t \; $ in probability as $n \rightarrow \infty$, so there is a subsequence $(n_k)$ such that \\
$\displaystyle A_{n_k,2} \rightarrow <S^H(N,a)>_t \; a.s.$ as $k \rightarrow \infty$ and therefore $\displaystyle \sup_k A_{n_k,2} < \infty \; a.s.$  Let
$$v(x) = \Bigg( \sup_k \sum_{j=1}^{n_k} \Big[ x\Big(\frac{jt}{n_k}\Big) - x\Big(\frac{(j-1)t}{n_k}\Big)  \Big]^2 \Bigg)^{1/2}, \; x \in C([0,t]),$$
$v$ is a measurable seminorm on $C([0,t])$ such that $v(S^H(N,a)) < \infty \; a.s.$ Then by Fernique's theorem (cf. \cite{FE}) $\mathbb{E} (v(S^H(N,a)))^q < \infty$ for all $q > 0$.  Pick $p > 2$ such that $\displaystyle H_{i_0} < \frac{1}{p} < \frac{1}{2}$. Then
\[ \begin{array}{rcl}
&& \displaystyle {\Bbb E} (A_{n_k,p})  \le {\Bbb E} \Bigg( \max _{j \le n} \Big| S^H_{j/n_k}(N,a)-S^H_{(j-1)/n_k}(N,a) \Big| ^{p-2} v(S^H(N,a))^2 \Bigg) \\
\noalign{\vskip 2mm}
&\le & \displaystyle \Bigg[ {\Bbb E} \Bigg( \max _{j \le n} \Big| S^H_{j/n_k}(N,a)-S^H_{(j-1)/n_k}(N,a) \Big| ^{2(p-2)} \Bigg) \Bigg] ^{1/2} \Big[ {\Bbb E} \Big( v(S^H(N,a))^4 \Big) \Big] ^{1/2} . \\
\end{array} \]
Since $p > 2$, the last expression tends to $0$ as $k \rightarrow \infty$; in fact, by continuity of $S^H(N,a)$, Lebesgue's theorem and the fact that
$${\Bbb E} \Big( \sup_{s \le t} \mid S^H_s(N,a) \mid ^{2(p-2)} \Big)  < \infty $$
( again by Fernique's theorem with the seminorm $\displaystyle x \longmapsto \sup_{s \le t} \mid x(s) \mid $ on $C([0,t])$),  we have
$${\Bbb E} \Bigg( \max _{j \le n} \Big| S^H_{j/n_k}(N,a)-S^H_{(j-1)/n_k}(N,a) \Big| ^{2(p-2)} \Bigg)  \rightarrow 0 \; as \; k \rightarrow \infty .$$
But, by (\ref{eq:15}) and since  the increments of $S^H(N,a)$ are Gaussian, it is easy to get
\begin{equation}
a n^{1- pH_{i_0}} \le E \Big( A_{n,p} \Big) 
\end{equation}
where $a$ is a positive constant, which depends on $\displaystyle a_{i_0}, H_{i_0}, p$ and $t$. Since $1 - pH_{i_0} > 0$, we deduce that $\displaystyle {\Bbb E} (A_{n,p}) $ tends to $\infty $ as $\displaystyle n \rightarrow \infty$. Hence it is not possible that $<S^H(N,a)>_t$ be finite a.s. and by the $0-1$ law  $<S^H(N,a)>_t = \infty \; a.s.$\\

$2.$ To check the first assertion, let us consider a sequence 
$$  \Big\{ \tau _n: 0 = t_0 < t_1 < ... < t_n = t,
 \Big\} _{ n \in {\Bbb N}\setminus \{ 0 \} } $$
 of finite partitions
of $\displaystyle [0,t]$  such that its mesh $\displaystyle \mid \tau _n \mid =
 \max_{i=1}^{n} \mid t_{i} - t_{i-1} \mid
\stackrel{ n \rightarrow  \infty }{\longrightarrow} 0$. From Proposition \ref{p:5},  for every $\displaystyle \sigma, s \in [ 0, t] $; 
$\displaystyle s \le \sigma $ \; we have \;
$$ \sum_{i=1}^N  a_i ^2 \gamma _i  (\sigma -s)^{2H_i}   \le  E \Big( S_\sigma^H(N,a) - S_s^H(N,a) \Big) ^2 \le \sum_{i=1}^N  a_i ^2 \nu _i (\sigma-s)^{2H_i} $$ 
and consequently, 
\begin{equation} \label{eq:68}
 C_1  (\sigma -s)^{2H_{i_0}}    \le  E \Big( S_\sigma^H(N,a) - S_s^H(N,a) \Big) ^2 \le C_2  (\sigma -s)^{2H_{i_0}} 
\end{equation}
where 
$$C_1= a_{i_0} ^2 \gamma _{i_0}\; \hspace{2mm} {\rm and }\; \hspace{2mm}  C_2 = \sum_{i=1}^N a_i^2 \nu _i t^{2(H_i-H_{i_0})}.$$

For  any integer $n \ge 1$, denote
$$\Delta _t^{\tau _n} = \sum_{j=1}^n ( S^H_{t_j}(N,a) - S^H_{t_{j-1}}(N,a) ) ^2 .$$
From  equation (\ref{eq:68}), 
 we have
\begin{equation} \label{eq:69}
 \begin{array}{rcl}
 \displaystyle E (  \Delta _t^{\tau _n}   ) &\le & \displaystyle  C_2 \sum_{j=1}^n (t_j -t_{j-1})^{2H_{i_0}}\\
 \noalign{\vskip 3mm}
 &\le & \displaystyle  C_2 \mid \tau _n \mid ^{2H_{i_0} -1} \sum_{j=1}^n (t_j -t_{j-1}) \\
 \noalign{\vskip 3mm}
 &= & \displaystyle  C_2 \mid \tau _n \mid ^{2H_{i_0} -1} t. 
 \end{array}
 \end{equation}

Since $\displaystyle  \lim_{n \rightarrow \infty  } C_2 \mid \tau _n \mid ^{2H_{i_0} -1} t = 0$, \;
$\displaystyle \lim_{n \rightarrow \infty } E (\Delta _t^{\tau _n}) = 0 . $

Hence the sequence $(\Delta _t^{\tau _n})$ converges to $0$ in probability, which yields that 
$$<S^H(N,a) > _t = 0 \; a.s. $$

In order to get the second assertion in $2.$,  it suffices to follow the same procedure as that of the proof of $1.$ in Lemma \ref{L:20}, and to use Fernique's Theorem and equation (\ref{eq:15}).
\qed

As a consequence of lemma  \ref{L:20} we get:

\begin{cor}  \label{C:21}
If it exists $\displaystyle i \in \{ 1, ..., N \}; H_i < 1/2 $ and $\displaystyle a_i \neq 0$, \; or
\;  if $\displaystyle H_i  > 1/2 $, for every $\displaystyle i \in \{ 1, ..., N \}$, 
then the msfBm $S^H(N,a)$ is not a semimartingale.

\end{cor}

\proof In the case where it exists $\displaystyle i \in \{ 1, ..., N \}; H_i < 1/2 $ and $\displaystyle a_i \neq 0$, the corollary is a direct consequence  of the assertion $1.$ of Lemma \ref{L:20}. So let us check the second case. Suppose $S^H(N,a)$ is  a right-continuous semimartingale. Hence, $S^H(N,a)$ can be written in the form
\begin{equation} 
\label{eq:72}
S^H(N,a) = M_t^H + A_t^H
\end{equation}
where $\displaystyle M^H_0 = A^H_0 = 0,$ $M^H$ is an a.s. right-continuous local martingale with respect to 
$\displaystyle \mathcal{\overline{F}}^{S^H(N,a)}$ and $A^H$ an a.s. right-continuous, $\displaystyle \mathcal{\overline{F}}^{S^H}-$ adapted finite variation process.
It follows, by  Theorem $II.22$ of \cite{Prott} and Lemma \ref{L:20}, that, for every $t \in [0;1]$,
$$0= < S^H(N,a), S^H(N,a) >_t = < M^H, M^H >_t.$$
By Theorem $II.27$ of \cite{Prott}, $M^H$ is itself  a zero process; and hence $S^H(N,a) = A^H$ has finite variation. This contradicts assertion $2.$ of Lemma \ref{L:20}. \qed

In the following lemma, we  treat the case where it exists $\displaystyle k_0 \in \{ 1, .., N \}; H_{k_0} = 1/2$ such that  $\displaystyle a_{k_0} \neq 0$ and for every 
$\displaystyle i \in \{ 1,..., N \} \setminus \{ k_0 \} $, $\displaystyle H_i  \in \big\{ 1/2 \big\} \cup ]3/4, 1[$. \\

\begin{lem} \label{L:24} 
For every $T > 0$, if it exists $\displaystyle k_0 \in \{ 1, .., N \}; H_{k_0} = 1/2$; $\displaystyle a_{k_0} \neq 0$ and for every 
$\displaystyle i \in \{ 1,..., N \} \setminus \{ k_0 \} $, $\displaystyle H_i  \in \Big\{ 1/2 \Big\} \cup \Big] 3/4 , 1 \Big[$, the process 
$$S^H(N,a) = \{ S_t^H(N,a), t \in [0,T] \} $$
is, in its own filtration,  a semimartingale equivalent in law with
$\displaystyle a_{k_0} $  times a Brownian motion.
\end{lem}

 \proof 
 The process $S^H(N,a) $ can be written \;
 $\displaystyle S^H(N,a) = a_{k_0} \Big( \xi^{1/2} + X_t \Big) $
 where
 $$X_t = \sum_{\stackrel{i=1}{i \neq k_0}}^N \frac{a_i}{a_{k_0}} \xi^{H_i}(t).$$

 The  process $\displaystyle X_t  $  is  Gaussian and its  covariance function 
$$R(s,t) = \sum_{\stackrel{i=1}{i \neq k_0}}^{N} \frac{a_i^2}{a_{k_0^2}}  \Bigg(   t^{2H_i} + s^{2H_i} - 1/2 \Big[ \mid s+t\mid ^{2H_i} + \mid t -s \mid ^{2H_i} \Big]
\Bigg) ,$$ 
is positive definite, twice continuously differentiable on $\displaystyle [0,T]^2 \setminus \{ (s,t); t = s \}$ and
for every $\displaystyle (s,t) \in [0,T]^2 \setminus \{ (s,t); t = s \}$,

$$ \frac{\partial ^2R(s,t) }{\partial s \partial t } = \sum_{\stackrel{i=1}{i \neq k_0}}^{N} \frac{a_i^2}{a_{k_0^2}} H_i(2H_i-1) \Big[ \mid t-s \mid ^{2H_i-2} - \mid s+t \mid ^{2H_i-2} \Big] .$$

So, because   $\displaystyle H_i \in \Big\{ 1/2 \Big\} \cup \Big] 3/4; 1 \Big[, $   for every $\displaystyle i \in \{1 , ..., N \} \setminus \{ k_0 \}$,   it's easy to check that $\displaystyle \frac{\partial ^2R}{\partial s \partial t } \in L^2([0,T]^2)$.\\

 On the other hand, $\displaystyle \xi^{1/2}$ is a Brownian motion independent of 
 $ \displaystyle X_t  $. 
 Consequently, according to \cite{Ba}, 
 the process 
 $$ \{ \xi_t^{1/2} + X_t, t \in [0,T] \}$$
 is, in its own filtration,  a semimartingale equivalent in law to a
 Brownian motion and the lemma is proved.
\qed

 We finish this section by the study of the last case. That is, the case where it exists $\displaystyle k_0 \in \{ 1, ..., N \}; H_{k_0} = 1/2$,  $H_i \ge 1/2$ for every $i \in \{ 1, ..., N \} \setminus \{ k_0 \}$, and it exists $i \neq k_0$; $ H_i \in \Big] 1/2, 3/4 \Big]$. We study a more general case,   where, for every $\displaystyle k \in \{ 1, ..., N \}, H_k \ge 1/2$ and   it exists $\displaystyle i \in \{ 1, ..., N \} ;  H_i \in ]1/2,3/4]$.  Let us first recall the definition of a quasimartingale. \\

 {\bf Definition:} 
 A stochastic process $(X_t)_{t \ge 0}$ is a quasimartingale if for every $\displaystyle T > 0, X_t \in L^1$ for all $ t \in [0,T],$ and 
 $$\sup_{\tau } \sum_{j=0}^{k-1} \left\| E \Big( X_{t_{j+1}} -X_{t_j} | {\cal F}_{t_j}^X \Big) \right\| _1 < \infty ,$$
 where $\tau $ is the set of all finite partitions
 $$
 0= t_0 < t_1 < ... < t_n = T \hspace{1mm} of \hspace{1mm} [0,T].
 $$

 \bigskip
 In the following key lemma we specify the relation between quasimartingale and semimartingale in the case of our process $S^H(N,a)$.

 \begin{lem}
 \label{L:25}
 If $ S^H(N,a)$ is not a quasimartingale, it is not a semimartingale, with respect to its own filtration.
 \end{lem}

 \proof The proof is similar to that of Lemma  $4.2$ in \cite{Cherid}; in fact,  only Gaussianity was used there.\qed

\begin{lem} \label{L:26} 
If, for all $\displaystyle k \in \{ 1, ..., N \}$,  $\displaystyle H_k \ge 1/2$, and it exists $\displaystyle k_0 \in \{ 1, ..., N \}$;  $\displaystyle   1/2 < H_{k_0} \le 3/4$ and $a_{k_0} \neq 0$, then $\displaystyle S^H(N,a)$ is not a quasimartingale, and in particular, it's not a semimartingale.
\end{lem}

 \proof 
The proof is an easy extention of the proof of Lemma $25$ in \cite{MZCE}  to our general case.
\qed

{\bf Author information}\\

 Mounir Zili, Department of Mathematics, 
Faculty of Sciences of Monastir, 
Avenue de l'environnement, $5019$ Monastir, Tunisia \\
Mounir.Zili@fsm.rnu.tn\\

\end{document}